\documentclass[11pt, english, a4paper]{amsart}
\usepackage[T1]{fontenc}
\usepackage[english]{babel}
\usepackage{amsmath,enumerate, amsfonts, amssymb, amsthm, graphicx}
\usepackage[all]{xy}

\newtheorem{defi}{Definition}[section]
\newtheorem{theo}[defi]{Theorem}
\newtheorem{prop}[defi]{Proposition}
\newtheorem{enonce}[defi]{enonce}
\newtheorem{rema}[defi]{Remark}
\newtheorem{lemm}[defi]{Lemma}
\newtheorem{coro}[defi]{Corollary}
\newtheorem{theorem}{Theorem}
\newtheorem*{lemm*}{Lemma}

\title{Irregularity of an analogue of the Gauss-Manin systems}
\author{C\'eline Roucairol}
\address{Laboratoire J. A. Dieudonn\'e, Universit\'e de Nice\\
parc Valrose, 06108 Nice cedex 02}
\email{roucair@math.unice.fr}

\begin{document}

\begin{abstract}
\noindent 
In the $\mathcal{D}$-modules theory, Gauss-Manin systems are defined by the direct image of the structure sheaf $\mathcal{O}$ by a morphism. A major theorem says that these systems have only regular singularities. This paper examines the irregularity of an analogue of the Gauss-Manin systems. It consists in the direct image complex of a $\mathcal{D}$-module twisted by the exponential of a polynomial $g$ by another polynomial $f$, $f_+(\mathcal{O}e^g)$, where $f$ and $g$ are two polynomials in two variables. The analogue of the Gauss-Manin systems can have irregular singularities (at finite distance and at infinity). We express an invariant associated with the irregularity of these systems at $c\in \mathbb{P}^1$ by the geometry of the map $(f,g)$. 
\end{abstract}

\maketitle
\tableofcontents
\section{Introduction}
\subsection{}
We denote by $\mathcal{O}_{\mathbb{C}^n}$ the sheaf of regular functions on 
$\mathbb{C}^n$ and by $\mathcal{D}_{\mathbb{C}^n}$ the sheaf of algebraic 
differential operators on $\mathbb{C}^n$. 

If $f:\mathbb{C}^n\to\mathbb{C}$ is a polynomial, we define the Gauss-Manin connection as the extension of the flat bundle $\underset{t\in\mathbb{C}\setminus\Sigma}{\cup}H^{k+n-1}(f^{-1}(t)^{an},\mathbb{C})$, where $\Sigma\subset\mathbb{C}$ is a finite subset such that $f:f^{-1}(\mathbb{C}\setminus\Sigma)\to\mathbb{C}\setminus\Sigma$ is a locally trivial fibration. A major theorem says that these connections are regular. In the $\mathcal{D}$-module theory, we study this connection with the help
of a complex of $\mathcal{D}_{\mathbb{C}}$-modules, it being the
direct image complex $f_+(\mathcal{O}_{\mathbb{C}^n})$. Their cohomology modules are called Gauss-Manin systems. They are holonomic and regular.

Now, let $g:\mathbb{C}^n\to\mathbb{C}$ be another polynomial. We denote by
$\mathcal{O}_{\mathbb{C}^n}e^g$ the
$\mathcal{D}_{\mathbb{C}^n}$-module obtained from
$\mathcal{O}_{\mathbb{C}^n}$ by twisting by $e^g$. We are interested in an
analogue of the Gauss-Manin systems, it being the direct image
complex $f_+(\mathcal{O}_{\mathbb{C}^n}e^g)$. 

In \cite{Ma}, F. Maaref calculates the generic fibre of the sheaf of horizontal analytic sections of the systems $\mathcal{H}^k(f_+(\mathcal{O}_{\mathbb{C}^n}e^g))$. It consists in a 
relative version of a result of C. Sabbah in \cite{Sa}. Indeed,
the generic fiber of the sheaf of horizontal analytic sections of $\mathcal{H}^k(f_+(\mathcal{O}_{\mathbb{C}^n}e^g))$ is canonically
isomorphic to the 
cohomology group with closed support
$H^{k+n-1}_{\Phi_t}(f^{-1}(t)^{an},\mathbb{C})$, where 
$\Phi_t$ is a family of closed subsets of $f^{-1}(t)$, on which
$e^{-g}$ is rapidly decreasing. More precisely, this family is defined
as follow. Let $\pi:\widetilde{\mathbb{P}^1}\to\mathbb{P}^1$ be the
oriented real blow-up of $\mathbb{P}^1$ at
infinity. $\widetilde{\mathbb{P}^1}$ is diffeomorphic
to $\mathbb{C}\cup S^1$, where $S^1$ is the circle of directions at
infinity. $A$ is in $\Phi_t$ if $A$ is a closed subset
of $f^{-1}(t)$ and the closure of $g(A)$ in $\mathbb{C}\cup S^1$
intersects $S^1$ in $]-\frac{\pi}{2},\frac{\pi}{2}[$.

This isomorphism can be better understood using relative
cohomology group. F. Maaref shows that for all $t\notin\Sigma$ and for
all $\rho$, such that $Re(-\rho)$ is sufficiently large, the fibre at
$t$ of the sheaf of horizontal analytic sections of
$\mathcal{H}^k(f_+(\mathcal{O}_{\mathbb{C}^n}e^g))$ is isomorphic to
the relative cohomology group $H^{k+n-1}(f^{-1}(t)^{an}, (f^{-1}(t)\cap
g^{-1}(\rho))^{an},\mathbb{C})$. 

Finally, he proves the quasi-unipotence of the corresponding
local mono\-dromy.

\subsection{}
The Gauss-Manin systems have only regular singularities. In our case, the complex $f_+(\mathcal{O}_{\mathbb{C}^n}e^g)$ can
have irregular singularities. The aim of this paper is to characterize this irregularity in terms of the geometry of the map $(f,g)$, when $f$ and $g$ are two polynomials in two variables. 

In $f$ and $g$ are algebraically independant, we will prove that the complex $f_+(\mathcal{O}_{\mathbb{C}^2}e^g)$ is essentially concentrated in degree zero. Then, we can associate to this complex a system of differential equations in one variable. We want to calculate the irregularity number of this system at a point at finite distance and at infinity. 

Let $\mathbb{X}$ be a smooth projective compactification of $\mathbb{C}^2$ such that there exists $F,G:\mathbb{X}\to\mathbb{P}^1$, two meromorphic maps, which extend $f$ and $g$. Let us denote by $D$ the divisor $\mathbb{X}\setminus\mathbb{C}^2$. In the following, we identify $\mathbb{P}^1$ with $\mathbb{C}\cup\{\infty\}$. 

Let $\Gamma$ be the critical locus of $(F,G)$. We denote by $\Delta_1$ the cycle in $\mathbb{P}^1\times\mathbb{P}^1$ which is the closure in
$\mathbb{P}^1\times\mathbb{P}^1$ of $(F,G)(\Gamma)\cap(\mathbb{C}^2\setminus\{c\}\times\mathbb{C})$ where the image is counted with multiplicity and by
$\Delta_2$ the cycle in $\mathbb{P}^1\times\mathbb{P}^1$ which is the closure in $\mathbb{P}^1\times\mathbb{P}^1$ of
$(F,G)(D)\cap(\mathbb{C}^2\setminus\{c\}\times\mathbb{C})$ where the image is counted with multiplicity.

For all $c\in\mathbb{P}^1$, the germs at $(c,\infty)$ of the support of $\Delta_1$ and $\Delta_2$ are some germs of curves or are empty. Then, we denote by $I_{(c,\infty)}(\Delta_i,\mathbb{P}^1\times\{\infty\})$ the intersection number of the cycles $\Delta_i$ and $\mathbb{P}^1\times\{\infty\}$. If the germ at $(c,\infty)$ of $\Delta_i$ is empty, this number is equal to $0$. 

\begin{theorem}
Let $f,g\in\mathbb{C}[x,y]$ be algebraically independant. Let $c\in\mathbb{P}^1$.\\
Then, the irregularity number at $c$ of the system $\mathcal{H}^0(f_+(\mathcal{O}_{\mathbb{C}^2}e^g))$ is equal to $I_{(c,\infty)}(\Delta_1,\mathbb{P}^1\times\{\infty\})+I_{(c,\infty)}(\Delta_2,\mathbb{P}^1\times\{\infty\})$.
\end{theorem}

When $c\in\mathbb{C}$, we can prove that the germ at $(c,\infty)$ of $\Delta_2$ is empty. Moreover, the germ at $(c,\infty)$ of $\Delta_1$ coincide with the one of the closure in $\mathbb{P}^1\times\mathbb{P}^1$ of $(f,g)(\tilde{\Gamma})\setminus\{c\}\times\mathbb{C}$, where $\tilde{\Gamma}$ is the critical locus of $(f,g)$.

\subsection{}
In general, we do not know how to calculate directly the irregularity number of a system associated with $f_+(\mathcal{O}_{\mathbb{C}^2}e^g)$. The notion of irregularity complex along an hypersurface defined by Z. Mebkhout (see \cite{Me1} and \cite{Me2}) is the appropriate tool to express the irregularity of $f_+(\mathcal{O}_{\mathbb{C}^2}e^g)$ (see paragraph \ref{irr}). Indeed, this irregularity complex along an hypersurface is a generalization of the irregularity number at a point of a system of differential equations in one variable. Moreover, Z. Mebkhout proves a theorem of commutation between the direct image functor and the irregularity functor (see theorem \ref{commute}). Then, the irregularity number at $c\in\mathbb{P}^1$ of the system of differential equations associated with $f_+(\mathcal{O}_{\mathbb{C}^2}e^g)$ can be expressed with the help of an irregularity complex of a $\mathcal{D}$-module in two variables along a curve.

In the general case where $f$ and $g$ are not necessarily algebraically independant, the complex $f_+(\mathcal{O}_{\mathbb{C}^2}e^g)$ is not necessarily concentrated in degree $0$. Then, we want to calculate the alternative sum of the irregularity number at $c\in\mathbb{P}^1$ of the systems $\mathcal{H}^k(f_+(\mathcal{O}_{\mathbb{C}^2}e^g))$. This irregularity number $IR_c$ is equal to the Euler characteristic of a complex of vector spaces over $\mathbb{C}$, it being the irregularity complex of $f_+(\mathcal{O}_{\mathbb{C}^2}e^g)$ at $c\in\mathbb{P}^1$. When $f$ and $g$ are algebraically independant, this number coincide with the irregularity number of the system $\mathcal{H}^0(f_+(\mathcal{O}_{\mathbb{C}^2}e^g))$. Then, we can prove that the irregularity number $IR_c$ is equal to $-\chi(\mathbb{R}\Gamma(F^{-1}(c)\cap G^{-1}(\infty), IR_{F^{-1}(c)}(\mathcal{O}_{\mathbb{X}}[\ast D]e^G)))$, where $IR_{F^{-1}(c)}(\mathcal{O}_{\mathbb{X}}[\ast D]e^G)$ is the irregularity complex of $\mathcal{O}_{\mathbb{X}}[\ast D]e^G$ along $F^{-1}(c)$. 

Then, according to a result of C. Sabbah in \cite{Sa}, we know that, for $x\in F^{-1}(c)\cap G^{-1}(\infty)$, the Euler characteristic of $(IR_{F^{-1}(c)}(\mathcal{O}_{\mathbb{X}}[\ast D]e^G))_x$ is equal to the Euler characteristic of the fiber $f^{-1}(D^{\ast}(c,\eta))\cap g^{-1}(\rho)\cap B(x,\epsilon)$, where $\epsilon$ and $\eta$ are small enough and $|\rho|$ is big enough. This result is stated in theorem \ref{Sab} paragraph \ref{exp} in terms of complex of nearby cycles. 

Then, we have to globalize the situation (see paragraph \ref{topo}). First of all, we prove that for $\eta$ small enough and $R$ big enough, $g:f^{-1}(D^{\ast}(c,\eta))\cap g^{-1}(\{|\rho|>R\})\to \{|\rho|>R\}$ is a locally trivial fibration. Then, the irregularity number $IR_c$ is equal to the opposite of the Euler characteristic of its fiber $f^{-1}(D^{\ast}(c,\eta))\cap g^{-1}(\rho)$. This result hold in the general case where $f$ and $g$ are not necessarily algebraically independant.

Then, we have to study the topology of this fiber. We have to distinguished the case where $f$ and $g$ are algebraically independant (see paragraph \ref{ind}) and the one where they are algebraically dependant (see paragraph \ref{dep}).

\section{Irregularity complex along an hypersurface}\label{irr}

We will use the definition of regularity given by Z. Mebkhout (see \cite{Me1} and \cite{Me2}). First of all, we recall the definition of irregularity complex of analytic $\mathcal{D}$-modules. Then, we define the notion of irregularity complex for algebraic $\mathcal{D}$-modules. Here, we have to take into account the behaviour of these modules at infinity. Moreover, we state major theorems on irregularity: the positivity theorem, the stability of the category of complex of regular holonomic $\mathcal{D}$-modules (analytic) by direct image by a proper map and the comparison theorem of Grothendieck.

\subsection{The analytic case}
Let $X$ be a smooth analytic variety over $\mathbb{C}$. In this section, $\mathcal{D}_X$ denotes the sheaf of analytic differential operators on $X$.

Let $Z$ be an analytic closed subset of $X$. Denote by $i$ the canonical inclusion of $X\setminus Z$ in $X$. Let $\mathcal{M}^{\bullet}$ be a bounded complex of analytic $\mathcal{D}_X$-modules with holonomic cohomology. 

\begin{defi}\label{irran}
We define the irregularity complex of $\mathcal{M}^{\bullet}$ along $Z$ as the complex :
$$\begin{array}{lll}
IR_Z(\mathcal{M}^{\bullet})&:=&
R\Gamma_Z(DR(\mathcal{M}^{\bullet}[\ast Z]))[+1]\\
&:=&cone\big(DR(\mathcal{M}[\ast Z])\to Ri_{\ast}i^{-1}(DR(\mathcal{M}^{\bullet}[\ast Z]))\big).
\end{array}$$
\end{defi}

According to the constructibility theorem (c.f. \cite{K} and \cite{Me3}), this complex is a bounded complex of constructible sheaves on $X$ with support in $Z$. Then, we can define the covariant exact functor $IR_Z$ between the category of bounded complexes of $\mathcal{D}_X$-modules with holonomic cohomology and the category of bounded complexes of constructible sheaves on $X$ with support in $Z$.

\begin{defi}
$\mathcal{M}^{\bullet}$ is said to be regular if its irregularity complex along all hypersurfaces of $X$ is zero. 
\end{defi}

In one variable, the previous definition of regularity generalises the notion of regular singular point of a differential equation which is characterized by the annulation of the irregularity number (Fuchs theorem). Indeed, irregularity complex along an hypersurface generalizes irregularity number in the case of one variable. According to Z. Mebkhout (see \cite{Me1}, \cite{Me2}), the characteristic cycle of the irregularity complex of a holonomic $\mathcal{D}$-module along an hypersurface is positive.

\begin{theo}[positivity theorem]\label{posit}
If $Z$ is an hypersurface of $X$ and $\mathcal{M}$ is a holonomic $\mathcal{D}_X$-module, the complex $IR_Z(\mathcal{M})$ is perverse on $Z$. 
\end{theo}

The category of complexes of $\mathcal{D}$-modules with regular holonomic cohomology is stable by proper direct image. Let us state the theorem which proves this stability (see \cite{Me1} and Prop. 3.6-4 of \cite{Me2}). It will be a major tool in this paper.

Let $\pi:X\to Y$ be a proper morphism of smooth analytic varieties over $\mathbb{C}$. Let $T$ be a hypersurface of $Y$. 

\begin{theo}\label{commute}
Let $\mathcal{M}^{\bullet}$ be a bounded complex of analytic $\mathcal{D}_X$-modules with holonomic cohomology. We have an isomorphism:
$$IR_T(\pi_+(\mathcal{M}^{\bullet}))[\dim Y]\simeq 
R\pi_{\ast}(IR_{\pi^{-1}(T)}(\mathcal{M}^{\bullet}))[\dim X].$$
\end{theo}

\subsection{The algebraic case}
Let $X$ be a smooth affine variety over $\mathbb{C}$. In this section, $\mathcal{D}_X$ denotes the sheaf of algebraic differential operators on $X$.

Denote by $j:X\to\mathbb{P}^n$ an immersion of $X$ in a projective space. Let $Z$ be a locally closed subvariety of $\mathbb{P}^n$ and $\mathcal{M}^{\bullet}$ a bounded complex of algebraic $\mathcal{D}_X$-modules with holonomic cohomology.

\begin{defi}\label{irral}
We define the irregularity complex of $\mathcal{M}^{\bullet}$ along $Z$ as the complex :
$$IR_Z(j_+(\mathcal{M}^{\bullet})):=IR_{Z^{an}}(j_+(\mathcal{M}^{\bullet})^{an}),$$
where $Z^{an}$ denotes the analytic variety associated with $Z$ and $j_+(\mathcal{M}^{\bullet})^{an}$ denotes the complex of analytic $\mathcal{D}$-modules associated with $j_+\mathcal{M}$. 
\end{defi}

\begin{defi}\label{regulier}
$\mathcal{M}^{\bullet}$ is said to be regular if its irregularity complex along all subvariety of $\mathbb{P}^n$ is zero. 
\end{defi}

This condition of regularity does not depend on the choice of the immersion $j$ (see proposition 9.0-4 in \cite{Me2}).

The definition \ref{regulier} of regular holonomic complex was motivated by the comparison theorem of Grothendieck (see \cite{Gro}) and the comparison theorem of Deligne (see \cite{Del}). As is shown in \cite{Me4}, the comparison theorem of Grothendieck is a consequence of the following theorem:

\begin{theo}\label{groth}
The structure sheaf $\mathcal{O}_X$ is regular in the sense of definition \ref{regulier}.
\end{theo}

Concerning to the stability of regularity under direct image, Theorem \ref{commute} allows to prove the following theorem (see Theorem 9.0-7 of \cite{Me2}):

\begin{theo}
The category of complexes of $\mathcal{D}$-modules with regular holonomic cohomology is stable under direct image.
\end{theo}

\begin{enonce}[remark]{Notation}
We denote by $IR_Z^k(j_+(\mathcal{M}^{\bullet}))$ the $k$-th space of cohomology of the complex
$IR_Z(j_+(\mathcal{M}^{\bullet}))$.
\end{enonce}

\begin{rema} We are interested in the irregularity of the direct image complex 
$f_+(\mathcal{O}_{\mathbb{C}^2}e^g)$. This is a complex of $\mathcal{D}_{\mathbb{C}}$-modules 
in one variable, with holonomic cohomology. According to the definition \ref{irral} of irregularity complex, we have 
to consider an immersion $j:\mathbb{C}\to\mathbb{P}^1$. Let $c\in\mathbb{P}^1$. 
We want to examine the complex $IR_{c}(j_+f_+(\mathcal{O}_{\mathbb{C}^n}e^g))$. 
As this complex has its support in $c$, we want to compute its Euler characteristic
$\chi(IR_{c}(j_+f_+(\mathcal{O}_{\mathbb{C}^n}e^g))_c)$. In the following, we will denote this number by $IR_c$. 
\end{rema}

\section{Regular holonomic $\mathcal{D}$-modules twisted by an exponential}\label{exp}

\subsection{Definitions}
Let $X$ be an algebraic variety over $\mathbb{C}$. We denote by $\mathcal{O}_X$ the sheaf of regular functions on $X$. 

We identify $\mathbb{P}^1$ to $\mathbb{C}\cup\{\infty\}$. Let $g:X\to\mathbb{P}^1$ be a meromorphic function on $X$. 
\begin{defi}
We define the $\mathcal{D}_X$-module $\mathcal{O}_X[\ast g^{-1}(\infty)]e^g$ as a $\mathcal{D}_X$-module which is isomorphic to $\mathcal{O}_{X}[\ast g^{-1}(\infty)]$ as $\mathcal{O}_X$-module; the action of $\xi$, vector field on an open subset of $X$, on a section $he^g$ of $\mathcal{O}_X[\ast g^{-1}(\infty)]e^g$ is defined by $\xi(he^g)=\xi(h)e^g+h\xi(g)e^g$. 
\end{defi}

Let $\mathcal{M}$ be a holonomic $\mathcal{D}_X$-module.

\begin{defi}
We define the $\mathcal{D}_X$-module $\mathcal{M}[\ast g^{-1}(\infty)]e^g$ as the $\mathcal{D}_X$-module $\mathcal{M}\otimes_{\mathcal{O}_X}\mathcal{O}_X[\ast g^{-1}(\infty)]e^g$.
\end{defi}

\begin{rema}
$\mathcal{O}_X[\ast g^{-1}(\infty)]e^g$ is the direct image by an open immersion of a vector bundle with integrable connection. Then, it is a holonomic $\mathcal{D}_X$-module as algebraic direct image of a holonomic $\mathcal{D}$-module. 

$\mathcal{M}[\ast g^{-1}(\infty)]e^g$ is a holonomic left $\mathcal{D}_X$-module as tensor product of two holonomic left $\mathcal{D}_X$-modules.
\end{rema}

We have analogous definitions in the analytic case. We just have to transpose in the analytic setting. 

\subsection{On irregularity of regular holonomic $\mathcal{D}$-modules twisted by an exponential}
Let $X$ be a complex analytic manifold and let $f,g:X\to\mathbb{C}$ be two
analytic functions. Assume that $\mathcal{M}$ is a regular holonomic
$\mathcal{D}_X$-module (analytic). Generally, $\mathcal{M}[\frac{1}{g}]e^{\frac{1}{g}}$ is an irregular $\mathcal{D}_X$-module. We want to relate the irregularity complex of this module along $f^{-1}(0)$ with some topological data. 

\begin{lemm}\label{Sab}
The complex
$IR_{f=0}(\mathcal{M}[\frac{1}{g}]e^{\frac{1}{g}})$ and the complex of
nearby cycles $\Psi_g(DR(\mathcal{M}[\frac{1}{f}]))$ have the same
characteristic function on $f^{-1}(0)\cap g^{-1}(0)$. 
\end{lemm}

\begin{proof}[Proof of Lemma \ref{Sab}] According to corollary 5.2 of \cite{Sa}, this lemma is true in the case where $f$ and $g$ are the same function. 

Assume that $f$ and $g$ are not equal. Then, using the case where the two functions are equal, we remark that it is sufficient to prove that the complex 
 $IR_{f=0}(\mathcal{M}[\frac{1}{g}]e^{\frac{1}{g}})$ and the complex 
$IR_{g=0}(\mathcal{M}[\frac{1}{fg}]e^{\frac{1}{g}})$ have the same 
characte\-ristic function on $f^{-1}(0)\cap g^{-1}(0)$.

\begin{itemize}
\item[$\bullet$] Let us first prove that $IR_{f=0}(\mathcal{M}[\frac{1}{g}]e^{\frac{1}{g}})
=R\Gamma_{f=0}(IR_{g=0}(\mathcal{M}[\frac{1}{fg}]e^{\frac{1}{g}}))$. 

Let $X^{\ast}$ denote $X\setminus g^{-1}(0)$ and $\eta$ be the inclusion
of $X^{\ast}$ in $X$. By definition, we have :
$$\begin{array}{lll}
IR_{g=0}(\mathcal{M}[\frac{1}{fg}]e^{\frac{1}{g}})&=&
cone(DR(\mathcal{M}[\frac{1}{fg}]e^{\frac{1}{g}})\to
R\eta_{\ast}\eta^{-1}DR(\mathcal{M}[\frac{1}{fg}]e^{\frac{1}{g}}))\\
&=&cone(DR(\mathcal{M}[\frac{1}{fg}]e^{\frac{1}{g}})\to
R\eta_{\ast}(DR(\mathcal{M}[\frac{1}{f}])_{|X^{\ast}})).
\end{array}$$

Now, consider the following diagram :
$$\xymatrix{X^{\ast}\ar@{^(->}[r]^{\eta}&X\\
X^{\ast}\setminus f^{-1}(0)\ar@{^(->}[r]^{\eta^{'}}\ar@{^(->}[u]_j&
X\setminus f^{-1}(0).\ar@{^(->}[u]_{j^{'}}}$$

Then, since $\mathcal{M}$ is regular, according to the definition \ref{regulier} of regularity, we have :
$$\begin{array}{lll}
R\eta_{\ast}(DR(\mathcal{M}[\frac{1}{f}])_{|X^{\ast}}))&=&
R\eta_{\ast}Rj_{\ast}(DR(\mathcal{M})_{|X^{\ast}\setminus f^{-1}(0)})\\
&=&Rj_{\ast}^{'}R\eta_{\ast}^{'}
(DR(\mathcal{M})_{|X^{\ast}\setminus f^{-1}(0)}).
\end{array}$$

As $R\Gamma_{f=0}Rj_{\ast}^{'}=0$, we obtain :
$$\begin{array}{lll}
R\Gamma_{f=0}(IR_{g=0}(\mathcal{M}[\frac{1}{fg}]e^{\frac{1}{g}}))
&=&R\Gamma_{f=0}(DR(\mathcal{M}[\frac{1}{fg}]e^{\frac{1}{g}}))[+1]\\
&=&IR_{f=0}(\mathcal{M}[\frac{1}{g}]e^{\frac{1}{g}}).
\end{array}$$

\item[$\bullet$] Then, we are led to show that the complex 
$R\Gamma_{f=0}(IR_{g=0}(\mathcal{M}[\frac{1}{fg}]e^{\frac{1}{g}}))$
and the complex $IR_{g=0}(\mathcal{M}[\frac{1}{fg}]e^{\frac{1}{g}})$
have the same characteristic function on $f^{-1}(0)\cap g^{-1}(0)$. 
Using the following distinguished triangle, 

\footnotesize
$$\xymatrix{&
Rj_{\ast}^{'}j^{'-1}(IR_{g=0}(\mathcal{M}[\frac{1}{fg}]e^{\frac{1}{g}}))
\ar[dl]_{[+1]}&\\
R\Gamma_{f=0}(IR_{g=0}(\mathcal{M}[\frac{1}{fg}]e^{\frac{1}{g}}))\ar[rr]&&
IR_{g=0}(\mathcal{M}[\frac{1}{fg}]e^{\frac{1}{g}}),\ar[ul]}$$
\normalsize
it is sufficient to show that the characteristic function of the complex
$Rj_{\ast}^{'}j^{'-1}(IR_{g=0}(\mathcal{M}[\frac{1}{fg}]e^{\frac{1}{g}}))$
is zero on $f^{-1}(0)\cap g^{-1}(0)$. Now, if $\mathcal{F}$ is a constructible sheaf on $X$ and $x\in f^{-1}(0)$, $\chi((Rj_{\ast}^{'}j^{'-1}\mathcal{F})_x)=\chi((D(j_!^{'}j^{'-1}D\mathcal{F}))_x)=\chi((j_!^{'}j^{'-1}D\mathcal{F})_x)=0$ ($D$ is the Verdier duality (see \cite{duality})). 
\end{itemize}
\end{proof}

\section{Topological interpretation of the irregularity of $f_+(\mathcal{O}_{\mathbb{C}^2}e^g)$}\label{topo}
\subsection{Notations}

Let $f,g:\mathbb{C}^2\to \mathbb{C}$ be two polynomials. 

Let $\mathbb{X}$ be a smooth projective compactification of $\mathbb{C}^2$ such that there exists $F,G:\mathbb{X}\to\mathbb{P}^1$ two meromorphic maps which extend $f$ and $g$. 

In view to construct $\mathbb{X}$, $F$ and $G$, we consider an immersion of $\mathbb{C}^2$ in $\mathbb{P}^2$ and we define a rational map $(\tilde{f},\tilde{g})$ on $\mathbb{P}^2$ which extends the map $(f,g)$. Then, after a finite number of blowing ups, we lift the indeterminacies of the rational map $(\tilde{f},\tilde{g})$.

In the following, we fixe such a compactification and use the following notations:
$$\xymatrix{\mathbb{C}^2\ar[r]^f\ar@{^(->}[d]^i&\mathbb{C}\ar@{^(->}[d]^j&&\mathbb{C}^2\ar[r]^g\ar@{^(->}[d]^i&\mathbb{C}\ar@{^(->}[d]^j\\
\mathbb{X}\ar[r]^F&\mathbb{P}^1&,&\mathbb{X}\ar[r]^G&\mathbb{P}^1.}$$

\subsection{Two fibration theorems}

\begin{lemm}\label{fibre1}
Let $c\in\mathbb{C}$. There exists $R>0$ big enough such that
$$g:g^{-1}(\{|\rho|>R\})\setminus (f^{-1}(c)\cap g^{-1}(\{|\rho|>R\}))\to \{|\rho|>R\}$$
is a locally trivial fibration.
\end{lemm}

\begin{proof}
Let $\mathcal{S}$ be an algebraic Whitney stratification of $\mathbb{X}$ such that $D$ and $F^{-1}(c)$ are union of strata. According to the Sard's theorem, there exists $U$, a dense Zariski open subset of $\mathbb{P}^1$, such that $G:G^{-1}(U)\to U$ is transverse to the Whitney stratification $\mathcal{S}^{'}$ of $G^{-1}(U)$ induced by $\mathcal{S}$. 

According to the first isotopy lemma of Thom-Mather, $G:G^{-1}(U)\to U$ is a locally trivial fibration with respect to $\mathcal{S}^{'}$. Then, we choose $R>0$ big enough such that $\{|\rho|>R\}\subset U$. 
\end{proof}

Using the inclusion $j:\mathbb{C}\to\mathbb{P}^1$, we identify $P^1$ to $\mathbb{C}\cup\{\infty\}$. If $c\in\mathbb{C}$, we denote by $D(c,\eta)$ the open disc in $\mathbb{C}$ centered at $c$ of radius $\eta$. Let $D(\infty,\eta)=\{z\in\mathbb{C}~|~|z|>\frac{1}{\eta}\}\cup\{\infty\}$. If $c\in\mathbb{P}^1$, we denote by $D^{\ast}(c,\eta)\subset\mathbb{C}$
the set $D(c,\eta)\setminus\{c\}$.

\begin{lemm}\label{fibre}
Let $c\in\mathbb{P}^1$. There exists $\eta$ small enough and $R$ big enough such that:
$$g:f^{-1}(D^{\ast}(c,\eta))\cap g^{-1}(\{|\rho|>R\})\to\{|\rho|>R\}$$
is a locally trivial fibration. 
\end{lemm}

\begin{proof}
According to the first isotopy lemma of Thom-Mather, we want to find a Whitney stratification $\mathcal{S}^{'}$ of $F^{-1}(\overline{D(c,\eta)})$ such that $D\cap F^{-1}(\overline{D(c,\eta)})$, $F^{-1}(c)$ and $F^{-1}(S(c,\eta))$ are union of strata and such that the morphism $G:F^{-1}(\overline{D(c,\eta)})\cap G^{-1}(\{|\rho|>R\})\to \{|\rho|>R\}$ is transverse to $\mathcal{S}^{'}$.

Let $\mathcal{S}$ be an algebraic Whitney stratification of $\mathbb{X}$ such that $F^{-1}(c)$ and $D$ are union of strata. For $\eta>0$, we denote by $\mathcal{T}$ the real analytic stratification $\{F^{-1}(S(c,\eta)), F^{-1}(D(c,\eta))\}$. For $\eta>0$ small enough, $\mathcal{S}$ and $\mathcal{T}$ are transverse. Let $\mathcal{S}^{'}$ be the real analytic stratification $\mathcal{S}\cap\mathcal{T}$. 

Now, let us prove that for $\eta$ small enough and $R$ big enough, the map $G:F^{-1}(\overline{D(c,\eta)})\cap G^{-1}(\{|\rho|>R\})\to \{|\rho|>R\}$ is transverse to $\mathcal{S}^{'}$.
\begin{itemize}
\item[$\bullet$] If $S^{'}=S\cap F^{-1}(S(c,\eta))$, for a $S\in\mathcal{S}$, we have to prove that for $\eta$ small enough and $R$ big enough, $G_{|S\cap F^{-1}(S(c,\eta))}$ is a submersion. It is sufficient to prove that for $\eta$ small enough and $R$ big enough, $F^{-1}(c^{'})$ is transverse to $G_{|S}^{-1}(\rho)$, where $|\rho|>R$ and $c^{'}\in S(c,\eta)$. 

Let $\Gamma_S$ be the critical locus of $(F,G)_{|S}$ and $\Delta_S=(F,G)(\Gamma_S)$ be the discriminant variety of $F_{|S}$ and $G_{|S}$. We denote by $\Delta_S^{'}$ the closure in $\mathbb{P}^1\times\mathbb{P}^1$ of $\Delta_S\cap\mathbb{C}^2$. In our case, the dimension of $\Delta_S^{'}$ is always less than $1$. 
Then we argue by the absurd. 
\item[$\bullet$] If $S^{'}=S\cap F^{-1}(D(c,\eta))$, for a $S\in\mathcal{S}$, as for $R$ big enough, the map $G:S\cap G^{-1}(\{|\rho|>R\})\to\{|\rho|>R\}$ is a submersion, the map $G:S^{'}\cap G^{-1}(\{|\rho|>R\})\to\{|\rho|>R\}$ is also a submersion.
\end{itemize}
\end{proof}

\subsection{Topological interpretation of the irregularity of $f_+(\mathcal{O}_{\mathbb{C}^2}e^g)$}
In this section, we describe the relation between the irregularity of the complex $f_+(\mathcal{O}_{\mathbb{C}^2}e^g)$ and the fibre $f^{-1}(D^{\ast}(c,\eta))\cap g^{-1}(\rho)$ given by the lemma \ref{fibre}.

\begin{theo}\label{theo1}
Let $c\in\mathbb{P}^1$. For $\eta$ small enough and $|\rho|$ big enough, 
$$IR_c=-\chi(f^{-1}(D^{\ast}(c,\eta))\cap g^{-1}(\rho)).$$
\end{theo}

This theorem can be proved in two steps. 
\begin{lemm}\label{L1}
$$IR_c=-\chi(\mathbb{R}\Gamma(F^{-1}(c)\cap G^{-1}(\infty),\psi_{\frac{1}{G}}(DR(\mathcal{O}_{\mathbb{X}}[\ast D\cup F^{-1}(c)])))).$$
\end{lemm}

\begin{lemm}\label{L2}
For $\eta$ small enough, $|\rho|$ big enough, $\chi(f^{-1}(D^{\ast}(c,\eta))\cap g^{-1}(\rho))$ is equal to $\chi(\mathbb{R}\Gamma(F^{-1}(c)\cap G^{-1}(\infty),\psi_{\frac{1}{G}}(DR(\mathcal{O}_{\mathbb{X}}[\ast D\cup F^{-1}(c)])))).$
\end{lemm}

\begin{proof}[Proof of lemma \ref{L1}]
This proof consists in applying the lemma \ref{Sab} on irregularity of regular holonomic $\mathcal{D}$-modules twisted by an exponential and in globalizing the situation. 

\begin{itemize}
\item[$\bullet$] First, we want to prove that
$$IR_c=-\chi(\mathbb{R}\Gamma(F^{-1}(c)\cap G^{-1}(\infty), IR_{F^{-1}(c)}(\mathcal{O}_{\mathbb{X}}[\ast D]e^G)^{an})).$$

According to the definition \ref{irral} and the theorem \ref{commute}, we have:
$$\begin{array}{lll}
IR_c(j_+f_+(\mathcal{O}_{\mathbb{C}^2}e^g))&=&
IR_c(F_+(\mathcal{O}_{\mathbb{X}}[\ast D]e^G))\\
&=&RF_{\ast}(IR_{F^{-1}(c)}(\mathcal{O}_{\mathbb{X}}[\ast D]e^G))[+1]
\end{array}$$

Then, $IR_c=-\chi(\mathbb{R}\Gamma(F^{-1}(c), IR_{F^{-1}(c)}(\mathcal{O}_{\mathbb{X}}[\ast D]e^G)^{an}))$. 
So we have to prove that the support of $IR_{F^{-1}(c)}(\mathcal{O}_{\mathbb{X}}[\ast D]e^G)^{an}$ is included in $F^{-1}(c)\cap G^{-1}(\infty)$. 

Let $x\notin G^{-1}(\infty)$. Then, $G$ is holomorphic in a neighbourhood of $x$ and $(\mathcal{O}_{\mathbb{X}}[\ast D]e^G)_x^{an}$ is isomorphic to $(\mathcal{O}_{\mathbb{X}}[\ast D])_x^{an}$. According to the theorem \ref{groth}, $\mathcal{O}_{\mathbb{X}}$ is regular. Then $\mathcal{O}_{\mathbb{X}}[\ast D]$ is also regular and $(IR_{F^{-1}(c)}(\mathcal{O}_{\mathbb{X}}[\ast D]e^G)^{an})_x=0$.  

\item[$\bullet$] According to the theorem \ref{Sab}, the complexes $IR_{F^{-1}(c)}(\mathcal{O}_{\mathbb{X}}[\ast D]e^G)^{an}$ and $\psi_{\frac{1}{G}}(DR(\mathcal{O}_{\mathbb{X}}[\ast D\cup F^{-1}(c)]))$ have the same characteristic function on $F^{-1}(c)\cap G^{-1}(\infty)$. We conclude using the following lemma:
\end{itemize}
\begin{lemm*}
Let $X$ be an algebraic variety over $\mathbb{C}$. Let $\mathcal{F}^{\bullet}_1$ and $\mathcal{F}^{\bullet}_2$
be two constructible complexes on $X$ which have the same characteristic 
function on $X$. Then,
$\chi(\mathbb{R}\Gamma(X,\mathcal{F}_1^{\bullet}))=
\chi(\mathbb{R}\Gamma(X,\mathcal{F}_2^{\bullet}))$.
\end{lemm*}

\begin{proof}
We argue by induction on the dimension of $X$. 
\begin{enumerate}
\item If the dimension of $X$ is $0$, the result is clear.
\item Assume that the lemma is true for all $X$ of dimension $<n$. Let $X$ be a complex algebraic manifold of dimension
$n$ and $Z$ be a closed complex algebraic submanifold of $X$, proper to $X$ ($dim
~Z<n$) such that $\mathcal{F}_{1|X\setminus Z}$ and 
$\mathcal{F}_{2|X\setminus Z}$ are some local systems $\mathcal{L}_1$ and $\mathcal{L}_2$ on $X\setminus Z$. For $i=1,2$, we have the
following distinguished triangle:
$$\xymatrix{\mathbb{R}\Gamma_Z(\mathcal{F}^{\bullet}_i)\ar[r]&
\mathcal{F}^{\bullet}_i\ar[r]&
\mathbb{R}j_{\ast}j^{-1}(\mathcal{F}^{\bullet}_i)\ar[r]^-{[+1]}&},$$
where $j$ is the inclusion of $X\setminus Z$ in $X$. Then:
$$\xymatrix{\mathbb{R}\Gamma(X,\mathbb{R}\Gamma_Z(\mathcal{F}^{\bullet}_i))
\ar[r]&\mathbb{R}\Gamma(X,\mathcal{F}^{\bullet}_i)\ar[r]&
\mathbb{R}\Gamma(X,\mathbb{R}j_{\ast}j^{-1}(\mathcal{F}^{\bullet}_i))
\ar[r]^-{[+1]}&}.$$
\begin{itemize}
\item[$\bullet$] $\chi(\mathbb{R}\Gamma(X,Rj_{\ast}j^{-1}(\mathcal{F}^{\bullet}_i)))=
\chi(\mathbb{R}\Gamma(X\setminus Z, \mathcal{L}_i))$.

As $X$ is an algebraic variety, $X\setminus Z$ is a finite union of connected open subsets $U_j$, $j=1,\ldots,k$, of $X$. Then, 
$$\chi(\mathbb{R}\Gamma(X,Rj_{\ast}j^{-1}(\mathcal{F}^{\bullet}_i)))=
\sum_{j=1}^k\chi(U_j)rk(\mathcal{L}_i).$$
As
$\mathcal{F}^{\bullet}_1$ and $\mathcal{F}^{\bullet}_2$ have the same
characteristic function, $rk(\mathcal{L}_1)=
rk(\mathcal{L}_2)$. Then,
$$\chi(\mathbb{R}\Gamma(X,Rj_{\ast}j^{-1}(\mathcal{F}^{\bullet}_1)))=
\chi(\mathbb{R}\Gamma(X,Rj_{\ast}j^{-1}(\mathcal{F}^{\bullet}_2))).$$
\item[$\bullet$] We have 
$\mathbb{R}\Gamma(X,\mathbb{R}\Gamma_Z(\mathcal{F}^{\bullet}_i))=
\mathbb{R}\Gamma(Z,\mathbb{R}\Gamma_Z(\mathcal{F}^{\bullet}_i))$. As at the end of the proof of lemma \ref{Sab}, using the Verdier duality (\cite{duality}), we can prove that $\mathbb{R}\Gamma(Z,\mathbb{R}\Gamma_Z(\mathcal{F}^{\bullet}_i))$ and $\mathbb{R}\Gamma(X,\mathcal{F}^{\bullet}_i)$ have the same characteristic function on $Z$. Then, $\mathbb{R}\Gamma(Z,\mathbb{R}\Gamma_Z(\mathcal{F}^{\bullet}_1))$ and $\mathbb{R}\Gamma(Z,\mathbb{R}\Gamma_Z(\mathcal{F}^{\bullet}_2))$ have the same characteristic function on $Z$. As $dim~Z<n$,
we apply the inductive hypothesis to obtain:
$$\chi(\mathbb{R}\Gamma(Y,\mathbb{R}\Gamma_Z(\mathcal{F}^{\bullet}_1)))=
\chi(\mathbb{R}\Gamma(Y,\mathbb{R}\Gamma_Z(\mathcal{F}^{\bullet}_2))).$$
\end{itemize}
Then,
$$\chi(\mathbb{R}\Gamma(Y,\mathcal{F}^{\bullet}_1))=
\chi(\mathbb{R}\Gamma(Y,\mathcal{F}^{\bullet}_2)).$$
\end{enumerate}
\end{proof}
\end{proof}

\begin{proof}[Proof of lemma \ref{L2}]
Denote by $\mathcal{F}^{\bullet}$ the complex $\psi_{\frac{1}{G}}(DR(\mathcal{O}_{\mathbb{X}}[\ast D\cup F^{-1}(c)]))$.
\begin{itemize}
\item[$\bullet$] Let $\eta_2$ small enough such that $G:G^{-1}(D^{\ast}(\infty,\eta_2))\to D^{\ast}(\infty,\eta_2)$ is a locally trivial fibration. We denote by $\widetilde{D^{\ast}(\infty,\eta_2)}$ the universal covering of $D^{\ast}(\infty,\eta_2)$. Let $(E,\pi,\widetilde{G})$ be the fiber product over $D^{\ast}(\infty,\eta_2)$ of $G^{-1}(D^{\ast}(\infty,\eta_2))$ and $\widetilde{D^{\ast}(\infty,\eta_2)}$. Then, we have the following diagram:

$$\xymatrix{G^{-1}(\infty)\ar@{^(->}[r]^-j&\mathbb{X}&G^{-1}(D^{\ast}(\infty,\eta_2))\ar@{_(->}[l]_-i\ar[d]_G&E\ar[l]_-{\pi}\ar[d]\\
&&D^{\ast}(\infty,\eta_2)&\widetilde{D^{\ast}(\infty,\eta_2)}\ar[l]}$$

By definition, $\mathcal{F}^{\bullet}=j^{-1}R(i\circ\pi)_{\ast}(i\circ\pi)^{-1}(DR(\mathcal{O}_{\mathbb{X}}[\ast D\cup F^{-1}(c)]))$. 

\item[$\bullet$] Let $\alpha:\mathbb{X}\setminus(F^{-1}(c)\cap D)\to \mathbb{X}$ open inclusion. As $\mathcal{O}_{\mathbb{X}}$ is regular in the sense of definition \ref{irral} (theorem \ref{groth}), 
$$DR(\mathcal{O}_{\mathbb{X}}[\ast D\cup F^{-1}(c)])=R\alpha_{\ast}\alpha^{-1}(\underline{\mathbb{C}}_{\mathbb{X}}).$$

\item[$\bullet$] Let $Z=F^{-1}(c)\cap G^{-1}(\infty)$, $Z_{\eta_1,\eta_2}=F^{-1}(D(c,\eta_1))\cap G^{-1}(D(\infty,\eta_2))$ and $Z_{\eta_1,\rho}=F^{-1}(D(c,\eta_1))\cap G^{-1}(\rho)$. 
$$\begin{array}{rl}
\mathbb{R}\Gamma(Z,\mathcal{F}^{\bullet})=&\underset{\eta_1,\eta_2>0}{indlim}\mathbb{R}\Gamma(Z_{\eta_1,\eta_2},R(i\circ\pi)_{\ast}(i\circ\pi)^{-1}(R\alpha_{\ast}\alpha^{-1}(\mathbb{C}_{\mathbb{X}}))\\
=& \underset{\eta_1,\eta_2>0}{indlim}\mathbb{R}\Gamma((i\circ\pi)^{-1}(Z_{\eta_1,\eta_2}),(i\circ\pi)^{-1}(R\alpha_{\ast}\alpha^{-1}(\mathbb{C}_{\mathbb{X}}))
\end{array}$$
Let $\Sigma$ be a Whitney stratification associated with the constructible sheaf $R\alpha_{\ast}\alpha^{-1}(\underline{\mathbb{C}}_{\mathbb{X}})$. Then, $F^{-1}(c)$ and $D$ are union of strata. According to the proof of lemma \ref{fibre}, for $\eta_1$ and $\eta_2$ small enough, $G:F^{-1}(D(c,\eta))\cap G^{-1}(D^{\ast}(\infty,\eta_2))\to D^{\ast}(\infty,\eta_2)$ is a locally trivial fibration with respect to $\Sigma$. Then, there exists a homotopy equivalence $p:(i\circ\pi)^{-1}(Z_{\eta_1,\eta_2})\to Z_{\eta_1,\rho}$ compatible with $\Sigma$. Thus, while adapting the proposition I.3-4 of \cite{Me3} to constructible sheaves, 
$$\begin{array}{rl}
\mathbb{R}\Gamma(Z,\mathcal{F}^{\bullet})=&\underset{\eta_1>0}{indlim}\mathbb{R}\Gamma(Z_{\eta_1,\rho},R\alpha_{\ast}\alpha^{-1}(\mathbb{C}_{\mathbb{X}}))\\
=&\underset{\eta_1>0}{indlim}\mathbb{R}\Gamma(\alpha^{-1}(Z_{\eta_1,\rho}),\alpha^{-1}(\mathbb{C}_{\mathbb{X}}))\\
=&\underset{\eta_1>0}{indlim}\mathbb{R}\Gamma(f^{-1}(D^{\ast}(c,\eta_1))\cap g^{-1}(\rho),\alpha^{-1}(\mathbb{C}_{\mathbb{X}}))\\
\end{array}$$
Then, $\chi(\mathbb{R}\Gamma(Z,\mathcal{F}^{\bullet}))=\chi(f^{-1}(D^{\ast}(c,\eta_1))\cap g^{-1}(\rho))$. 
\end{itemize}
\end{proof}

\section{When $f$ and $g$ are algebraically independant}\label{ind}

Let $f,g\in\mathbb{C}[x,y]$ be two polynomials which are algebraically independant. In this section, we will prove that the complex $f_+(\mathcal{O}_{\mathbb{C}^2}e^g)$ is essentially concentrated in degree $0$. Finally we obtain a formula for the irregularity number at $c\in\mathbb{P}^1$ in terms of some geometric data associated with $f$ and $g$.

\subsection{The complex $f_+(\mathcal{O}_{\mathbb{C}^2}e^g)$ is essentially concentrated in degree zero}

\begin{prop}
The complex $f_+(\mathcal{O}_{\mathbb{C}^2}e^g)$ is concentrated in degree zero except at a finite number of points. 
\end{prop}

\begin{proof}
\begin{itemize}
\item[$\bullet$] First of all, we recall the result of F. Maaref \cite{Ma} about the
gene\-ric fibre of the sheaf of horizontal analytic
sections of $\mathcal{H}^{k-1}(f_+(\mathcal{O}_{\mathbb{C}^2}e^g))$.

\begin{theo}\label{theoma}
There exists a finite subset $\Sigma$ of $\mathbb{C}$ such that for all $c\in\mathbb{C}\setminus\Sigma$ and
all $\rho\in\mathbb{C}$, such that $Re(-\rho)$ is big enough,
$$i_c^+\mathcal{H}^{k-1}(f_+(\mathcal{O}_{\mathbb{C}^2}e^g))\simeq H^k(f^{-1}(c),(f,g)^{-1}(c,\rho),\mathbb{C}),$$
where $i_c$ is the inclusion of $\{c\}$ in $\mathbb{C}$.
\end{theo}

\item[$\bullet$] For all $c,\rho\in\mathbb{C}$, we have the long exact sequence of relative cohomology:
\end{itemize}
\small
$$\xymatrix{0\ar[r]& H^0(f^{-1}(c),(f,g)^{-1}(c,\rho),\mathbb{C})\ar[r]&
H^0(f^{-1}(c),\mathbb{C})\ar[r]^-{\alpha}& H^0((f,g)^{-1}(c,\rho),\mathbb{C})\ar[dll] \\
&H^1(f^{-1}(c),(f,g)^{-1}(c,\rho),\mathbb{C})\ar[r]&
H^1(f^{-1}(c),\mathbb{C})\ar[r]^-{\beta}& H^1((f,g)^{-1}(c,\rho),\mathbb{C})\ar[dll]\\
&H^2(f^{-1}(c),(f,g)^{-1}(c,\rho)),\mathbb{C})\ar[r]& 0.}$$

\normalsize
\begin{itemize}
\item[] We want to prove that $H^k(f^{-1}(c),(f,g)^{-1}(c,\rho))=0$ for all $k\neq 1$. As $H^1((f,g)^{-1}(c,\rho),\mathbb{C})=0$, it is enough to prove that $\alpha$ is injective. Then, it is sufficient to prove that the fibre $g^{-1}(\rho)$ intersects all the connected components of $f^{-1}(c)$. 

Let $(F,G):\mathbb{X}\to\mathbb{P}^1\times\mathbb{P}^1$ be a compactification of $(f,g):\mathbb{C}^2\to\mathbb{C}^2$. As $(F,G):\mathbb{X}\to\mathbb{P}^1\times\mathbb{P}^1$ is proper, we know that its image is closed in $\mathbb{P}^1\times\mathbb{P}^1$. Moreover, as $f$ and $g$ are algebraically independant, $(F,G)$ is necessarily surjective. 

According to the theorem of Stein factorization (\cite{Ha} corollary 11.5 p. 280), there exists
$F^{'}:\mathbb{X}\to Y$, surjective morphism of projective varieties with connected fibres and
a finite morphism $\Psi:Y\to\mathbb{P}^1$, such that $F=\psi\circ F^{'}$. As $(F,G)$ is surjective, $(F^{'},G)$ is also surjective. Then, for all $(c,\rho)\in\mathbb{P}^1\times\mathbb{P}^1$, $G^{-1}(\rho)$ intersects all the connected components of $F^{-1}(c)$. 

Furthermore, there exists $\Sigma\subset\mathbb{C}$ finite subset such that for all $c\in\mathbb{C}\setminus\Sigma$, the fibre $F^{-1}(c)$ is the union of $f^{-1}(c)$ with a finite number of points. Then, for a such $c$, there exists $\Sigma_c\subset\mathbb{C}$ finite subset such that for all $\rho\in \mathbb{C}\setminus\Sigma_c$, $g^{-1}(\rho)$ intersects all the connected components of $f^{-1}(c)$. Then, for all $(c,\rho)\in\mathbb{C}^2$ except a finite number, the fibre $g^{-1}(\rho)$ intersects all the connected components of $f^{-1}(c)$. 

\item[$\bullet$] Then, according to the theorem \ref{theoma}, for all $c\in\mathbb{C}\setminus\Sigma$, for all
$k\neq 0$, $i_c^+(\mathcal{H}^{k}(f_+(\mathcal{O}_{\mathbb{C}^2}e^g))=0$. As
$\mathcal{H}^{k}(f_+(\mathcal{O}_{\mathbb{C}^2}e^g))$ is an integrable
connection except at a finite number of
points, we have that, for all $c\in\mathbb{C}$ except a finite number,
$\mathcal{H}^{k}(f_+(\mathcal{O}_{\mathbb{C}^2}e^g))_c=0$, if $k\neq 0$. Thus,
$f_+(\mathcal{O}_{\mathbb{C}^2}e^g)$ is essentially concentrated in degree $0$.

\end{itemize}
\end{proof}

\begin{coro}
For all $c\in\mathbb{P}^1$, the complex $IR_c(j_+f_+(\mathcal{O}_{\mathbb{C}^2}e^g))$ is concentrated in degree $0$.
\end{coro}

\begin{proof}
Let $c\in\mathbb{P}^1$. Denote by $\mathcal{M}^{\bullet}$ the complex $j_+f_+(\mathcal{O}_{\mathbb{C}^2}e^g)$. For $k\neq 0$, $\mathcal{H}^{k}(\mathcal{M}^{\bullet})$ has
punctual support. Let $\eta$ be small enough such that for all $k\neq 0$, $\mathcal{H}^k(\mathcal{M}^{\bullet})_{|D(c,\eta)^{\ast}}=0$. 
Then, $(\mathcal{M}^{\bullet}[\ast \{c\}])_{|D(c,\eta)}=(\mathcal{H}^0(\mathcal{M}^{\bullet})[\ast\{c\}])_{|D(c,\eta)}$.

Then, $IR_c(\mathcal{M}^{\bullet})=IR_c(\mathcal{H}^0(\mathcal{M}^{\bullet})[\ast\{c\}])$. 
Then, according to the positivity theorem \ref{posit}, this complex is just a vector space over $\mathbb{C}$. So, for all $k\neq 0$,
$\mathcal{H}^k(IR_c(j_+f_+(\mathcal{O}_{\mathbb{C}^2}e^g)))=0$.
\end{proof}

\begin{enonce}[remark]{Remark}
According to this corollary, the complex $IR_c(j_+f_+(\mathcal{O}_{\mathbb{C}^2}e^g))$ is entirely determined by its Euler characteristic $IR_c$. 
\end{enonce}

\subsection{Geometrical interpretation of the irregularity}

\begin{enonce}[remark]{Notation}
\begin{itemize}
\item[$\bullet$] Let $\Gamma$ be the critical locus of $(F,G)$. In the case where $f$ and $g$ are algebraically independant, this variety has dimension $1$.
\item[$\bullet$] Let $\Delta$ be the discriminant variety of $F$ and $G$. $\Delta$ is the image by $(F,G)$ of the curve $\Gamma$ (counted with multiplicity). $\Delta$ is an algebraic closed subset of $\mathbb{P}^1\times\mathbb{P}^1$. 
\item[$\bullet$] Denote by $\Delta_1$ the cycle in $\mathbb{P}^1\times\mathbb{P}^1$ which is the closure in $\mathbb{P}^1\times\mathbb{P}^1$ of $\Delta\cap(\mathbb{C}^2\setminus\{c\}\times\mathbb{C})$, where $\Delta$ is counted with multiplicity. 
\item[$\bullet$] Denote by $\Delta_2$ the cycle in $\mathbb{P}^1\times\mathbb{P}^1$ which is the closure in $\mathbb{P}^1\times\mathbb{P}^1$ of $(F,G)(D)\cap(\mathbb{C}^2\setminus\{c\}\times\mathbb{C})$, where the image is counted with multiplicity. 
\item[] (the supports of $\Delta_1$ and $\Delta_2$ are two algebraic closed subsets in $\mathbb{P}^1\times\mathbb{P}^1$. They are some union of curves and points.)
\item[$\bullet$] For all $c\in\mathbb{P}^1$, the germs at $(c,\infty)$ of the supports of $\Delta_1$ and $\Delta_2$ are some germs of curves or are empty. We denote by $I_{(c,\infty)}(\Delta_i,\mathbb{P}^1\times\{\infty\})$ the intersection number of the cycles $\Delta_i$ and $\mathbb{P}^1\times\{\infty\}$. If the germ at $(c,\infty)$ of $\Delta_i$ is empty, this number is equal to $0$. 
\end{itemize}
\end{enonce}

\begin{theo}
Let $f,g\in\mathbb{C}[x,y]$ be two polynomials algebraically independants. Let $c\in\mathbb{P}^1$. 
$$IR_c=I_{(c,\infty)}(\Delta_1,\mathbb{P}^1\times\{\infty\})+
I_{(c,\infty)}(\Delta_2,\mathbb{P}^1\times\{\infty\}).$$
\end{theo}

\begin{proof}
According to theorem \ref{theo1}, 
$$IR_c=-\chi(f^{-1}(D^{\ast}(c,\eta))\cap g^{-1}(\rho)).$$
We want to study the topology of the fibre $f^{-1}(D^{\ast}(c,\eta))\cap g^{-1}(\rho)$. 

For $|\rho|$ big enough, $G^{-1}(\rho)$ is smooth and cut transversally $D$. Then, 
$$\begin{array}{ll}
\chi(f^{-1}(D^{\ast}(c,\eta))\cap g^{-1}(\rho))=&
\chi(F^{-1}(D^{\ast}(c,\eta))\cap G^{-1}(\rho))\\
&-\chi(F^{-1}(D^{\ast}(c,\eta))\cap G^{-1}(\rho)\cap D).
\end{array}$$
\begin{itemize}
\item[$\bullet$] We want to prove that $\chi(F^{-1}(D^{\ast}(c,\eta))\cap G^{-1}(\rho))=I_{(c,\infty)}(\Delta_1,\mathbb{P}^1\times\{\infty\})$.

As $\Delta_1$ is a union of curves and points, for $\eta$ small enough and $|\rho|$ big enough, $\Delta_1\cap(D^{\ast}(c,\eta)\times\{\rho\})$ is a finite set $\{(c_1,\rho),\ldots,(c_r,\rho)\}$. Moreover, $F:F^{-1}(D^{\ast}(c,\eta))\cap G^{-1}(\rho)\to D^{\ast}(c,\eta)$ is a ramified covering. The ramified points are the points $P\in F^{-1}(D^{\ast}(c,\eta))\cap G^{-1}(\rho)\cap\Gamma$ and the ramification index at $P$ is $I_P(F^{-1}(F(P)),G^{-1}(\rho))=I_P(\Gamma, G^{-1}(\rho))+1$. Then,
$$\chi(F^{-1}(D^{\ast}(c,\eta))\cap G^{-1}(\rho))=
-\displaystyle\sum_{P\in F^{-1}(D^{\ast}(c,\eta))\cap G^{-1}(\rho)\cap\Gamma}I_P(\Gamma, G^{-1}(\rho)).$$
According to the projection formula for a proper map,
$$\begin{array}{lll}
\chi(F^{-1}(D^{\ast}(c,\eta))\cap G^{-1}(\rho))&=&
\displaystyle\sum_{i=1}^rI_{(c_i,\rho)}(\Delta,\mathbb{P}^1\times\{\rho\})\\
&=&I_{(c,\infty)}(\Delta_1,\mathbb{P}^1\times\{\infty\})
\end{array}$$
\item[$\bullet$] We want to prove that $\chi(F^{-1}(D^{\ast}(c,\eta))\cap G^{-1}(\rho)\cap D)=I_{(c,\infty)}(\Delta_2,\mathbb{P}^1\times\{\infty\})$.

As $G^{-1}(\rho)$ is smooth and cut transversally $D$,
$$\begin{array}{lll}
\chi(F^{-1}(D^{\ast}(c,\eta))\cap G^{-1}(\rho)\cap D)&=&
Card(F^{-1}(D^{\ast}(c,\eta))\cap G^{-1}(\rho)\cap D)\\
&=&\displaystyle\sum_{Q\in F^{-1}(D^{\ast}(c,\eta))\cap G^{-1}(\rho)\cap D}I_Q(D,G^{-1}(\rho))
\end{array}$$

According to the projection formula for a proper map,
$$\begin{array}{lll}
\chi(F^{-1}(D^{\ast}(c,\eta))\cap G^{-1}(\rho)\cap D)&=&
\displaystyle \sum_{(c^{'},\rho)\in \Delta_2\cap(D^{\ast}(c,\eta)\times\{\rho\})}I_{(c^{'},\rho)}((F,G)(D),\mathbb{P}^1\times\{\rho\})\\
&=&I_{(c,\infty)}(\Delta_2,\mathbb{P}^1\times\{\infty\})
\end{array}$$
\end{itemize}
\end{proof}

\begin{enonce}{Remark}
Denote by $\tilde{\Gamma}$ the critical locus of $(f,g)$. Let $\tilde{\Delta}$ be the closure of $(f,g)(\Gamma)\setminus (\{c\}\times\mathbb{C})$ in $\mathbb{P}^1\times\mathbb{P}^1$. 

If $c\in\mathbb{C}$, the germ at $(c,\infty)$ of $\Delta_1$ is the germ at $(c,\infty)$ of $\tilde{\Delta}$ and the one at $(c,\infty)$ of $\Delta_2$ is empty. Then, if $c\in\mathbb{C}$,
$$IR_c=I_{(c,\infty)}(\tilde{\Delta},\mathbb{P}^1\times\{\infty\}).$$
\end{enonce}

We deduce this remark from the following lemma:
\begin{lemm}
Let $c\in\mathbb{C}$. For $\eta>0$ small enough and $R>0$ big enough, 
$F^{-1}(D^{\ast}(c,\eta))\cap G^{-1}(\rho)$ does not intersect $D$, for $|\rho|>R$. 
\end{lemm}

\begin{proof}
We recall that $F,G:\mathbb{X}\to\mathbb{P}^1$ are constructed in the following way. We consider an immersion of $\mathbb{C}^2$ in $\mathbb{P}^2$. We construct a rational map $(\tilde{f},\tilde{g}):\mathbb{P}^2-->\mathbb{P}^1$ which extend $(f,g)$. On $\mathbb{P}^2\setminus\mathbb{C}^2$, $(\tilde{f},\tilde{g})$ takes the value $(\infty,\infty)$ or is not well defined. Then, we lift the indeterminacies of $(\tilde{f},\tilde{g})$ after a finite number of blowing ups. Then, according to \cite{LeWe}, we know that $F^{-1}(\infty)$ and $G^{-1}(\infty)$ are connected. As $F^{-1}(\infty)$ and $G^{-1}(\infty)$ have a non empty intersection, $F^{-1}(\infty)\cup G^{-1}(\infty)$ is connected.

Now let $Z$ be an irreducible component of $D$. We want to prove that for $\eta$ small enough and $R$ big enough, $F^{-1}(D^{\ast}(c,\eta))\cap G^{-1}(\rho)\cap Z=\emptyset$. 

If it is not true, we can construct a sequence $(x_n)_{n\in\mathbb{N}}$ such that $x_n\in Z$, $0<|F(x_n)-c|<\frac{1}{n}$ and $|G(x_n)|=\rho_n$, with $\lim_{n\to +\infty}\rho_n=\infty$. Then there exists a point $x\in F^{-1}(c)\cap G^{-1}(\infty)\cap Z$. 

As $Z\simeq \mathbb{P}^1$, $F_{|Z}$ and $G_{|Z}$ are necessarily surjective or constant. The existence of the sequence $(x_n)_{n\in\mathbb{N}}$ allows us to conclude that $F_{|Z}$ and $G_{|Z}$ are necessarily surjective. Then there exists another point $y\in Z\cap F^{-1}(\infty)$, $y\neq x$. 

This contradicts the facts that $F^{-1}(\infty)\cup G^{-1}(\infty)$ is connected and $F_{|Z}$ is surjective. 

\end{proof}

\section{When $f$ and $g$ are algebraically dependant}\label{dep}

Let $f,g\in\mathbb{C}[x,y]$ be two polynomials which are algebraically dependant. Then, the complex $f_+(\mathcal{O}_{\mathbb{C}^2}e^g)$ is not necessarily concentrated in degree $0$. However, we give a formula for the irregularity number $IR_c$ at $c\in\mathbb{P}^1$. 

Let $(F,G):\mathbb{X}\to\mathbb{P}^1\times\mathbb{P}^1$ be a compactification of the map $(f,g):\mathbb{C}^2\to\mathbb{C}^2$. As $f$ and $g$ are algebraically dependant, $\tilde{\Delta}=im~(F,G)$ is a closed subvariety of $\mathbb{P}^1\times\mathbb{P}^1$. 
\begin{itemize}
\item[$\bullet$] Let $\Delta$ be the cycle which is the closure in $\mathbb{P}^1\times\mathbb{P}^1$ of $\tilde{\Delta}\cap(\mathbb{C}^2\setminus\{c\}\times\mathbb{C})$. In a neighbourhood of $(c,\infty)$, $\Delta$ is a curve or is empty. We denote by $I_{(c,\infty)}(\Delta,\mathbb{P}^1\times\{\infty\})$ the intersection number of the cycles $\Delta$ and $\mathbb{P}^1\times\{\infty\}$. If the germ at $(c,\infty)$ of $\Delta$ is empty, this number is equal to $0$. 
\item[$\bullet$] Let $F$ be the generic fiber of $(f,g):\mathbb{C}^2\to im~(f,g)$. 
\end{itemize}

\begin{theo}
Let $f,g\in\mathbb{C}[x,y]$ be two polynomials algebraically dependants. Let $c\in\mathbb{P}^1$. 
$$IR_c=-\chi(F)*I_{(c,\infty)}(\Delta,\mathbb{P}^1\times\{\infty\}).$$
\end{theo}

\begin{proof}
According to the theorem \ref{theo1}, $IR_c=-\chi(f^{-1}(D^{\ast}(c,\eta))\cap g^{-1}(\rho))$, where $\eta$ is small enough and $|\rho|$ is big enough. 

As $\Delta$ is a curve or is empty in a neighbouhood of $(c,\infty)$, $\Delta\cap(D^{\ast}(c,\eta)\times\{\rho\})$ is a finite union of points $(c_1,\rho),\ldots,(c_r,\rho)$. 
Then, $f^{-1}(D^{\ast}(c,\eta))\cap g^{-1}(\rho)=\bigcup_{i=1}^r(f,g)^{-1}(c_i,\rho)$.

We conclude that $\chi(f^{-1}(D^{\ast}(c,\eta))\cap g^{-1}(\rho))=\chi(F)*I_{(c,\infty)}(\Delta,\mathbb{P}^1\times\{\infty\})$.
\end{proof}

\end{document}